\documentclass
[notitlepage,longbibliography,superscriptaddress,pre,nofootinbib]{revtex4-1}%
\usepackage{placeins}
\usepackage{pbox}
\usepackage{color}
\usepackage{relsize}
\usepackage{graphicx}
\usepackage[intlimits]{amsmath}
\usepackage{amsxtra, amssymb,fancyhdr, amsthm,latexsym}
\usepackage{cases}
\usepackage{verbatim}
\usepackage{float}
\usepackage[hidelinks]{hyperref}
\usepackage{textcomp}
\usepackage{amsfonts}
\usepackage{amssymb}%
\providecommand{\U}[1]{\protect\rule{.1in}{.1in}}
\newenvironment{absolutelynopagebreak}
{\par\nobreak\vfil\penalty0\vfilneg
\vtop\bgroup}
{\par\xdef\tpd{\the\prevdepth}\egroup
\prevdepth=\tpd}
\newtheorem{theorem}{Theorem}[section]

\newtheorem{mainres}[theorem]{Main Result}

\begin{document}
\begin{empty}
\title{It is better to be semi-regular when you have a low degree}
\author{T. Kolokolnikov}
\affiliation{Department of Mathematics and Statistics,
Dalhousie University Halifax,
Nova Scotia, B3H3J5, Canada}
\email{tkolokol@gmail.com}
\begin{abstract}
We study the algebraic connectivity for several classes of
random semi-regular graphs. For large random semi-regular bipartite graphs, we
explicitly compute both their algebraic connectivity and as well as the full
spectrum distribution. For an integer $d\in\left[  3,7\right]  $, we find
families of random semi-regular graphs that have higher algebraic connectivity
than a random $d$-regular graphs with the same number of vertices and edges.
On the other hand, we show that regular graphs beat semi-regular graphs when
$d\geq8.$ More generally, we study random semi-regular graphs whose average
degree is $d$, not necessary an integer. This provides a natural
generalization of a $d$-regular graph in the case of a non-integer $d.$ We
characterise their algebraic connectivity in terms of a root of a certain
6th-degree polynomial. Finally, we construct a small-world-type
network of average degree 2.5 with a relatively high algebraic connectivity.
We also propose some related open problems and conjectures.
\end{abstract}
\maketitle
\end{empty}

\section{Introduction}


Algebraic connectivity (AC; also called the spectral gap) of a graph is a
fundamental property that measures how fast information diffuses throughout
the graph \cite{fiedler1973algebraic, de2007old}. It corresponds to the second
smallest eigenvalue of the graph Laplacian matrix (the smallest eigenvalue is
always zero). In many applications, it is desirable to maximize the algebraic
connectivity (i.e. speed of diffusion)\ subject to certain constraints; this
and related problems have a long history \cite{ghosh2006growing,
wang2008algebraic, sydney2013optimizing, olfati2005ultrafast,
ghosh2006growing,kolokolnikov2015maximizing, ogiwara2015maximizing,
shahbaz2021algebraic}. For example in communications, the \textquotedblleft
cost\textquotedblright\ of a network increases with the number of links. It is
therefore desirable to have as few edges as possible, while at the same time
maximizing the algebraic connectivity. This leads to a natural question:
\emph{For a fixed average degree }$d,$\emph{ what is the graph that maximizes
algebraic connectivity as the number of vertices }$n\rightarrow\infty?$

\begin{figure}[tb]
\includegraphics[width=0.7\textwidth]{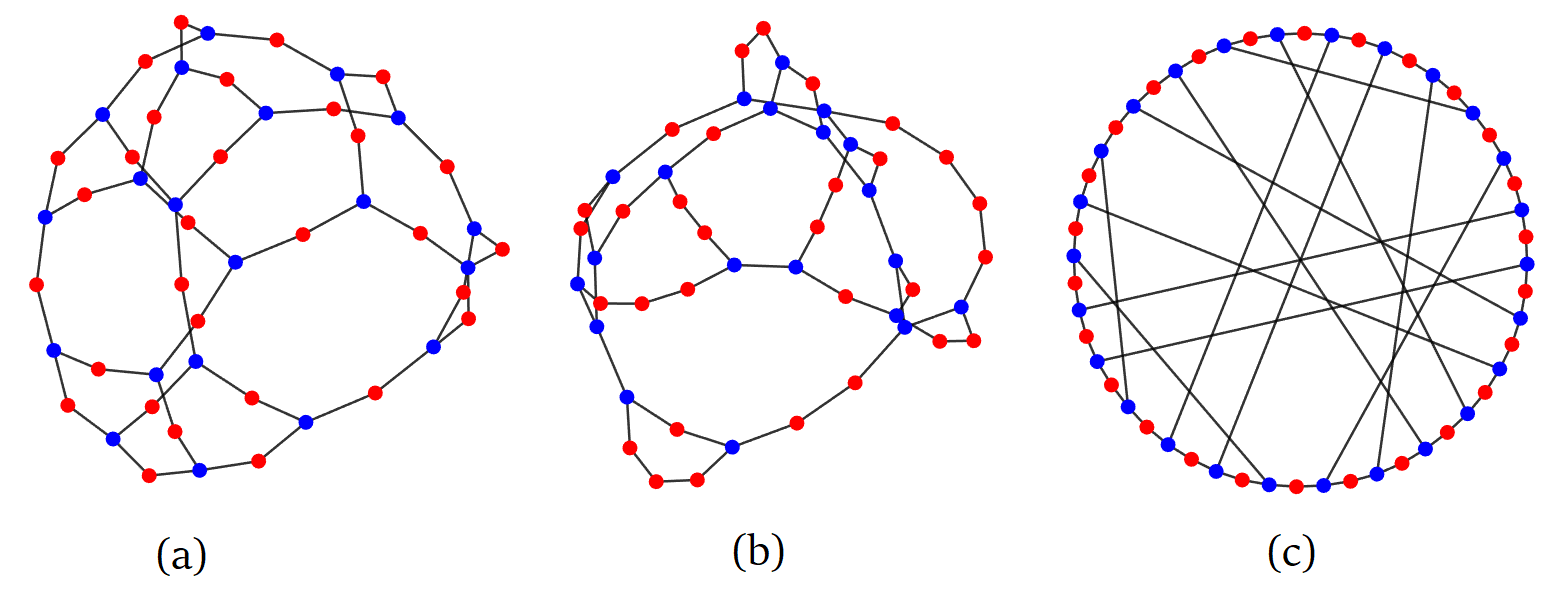}\caption{Types of graphs
considered in this work. (a)\ Random semi-regular bipartite graph with
$(d_{1},d_{2})=(2,3)$. (b) Random semi-regular graph with $(p,d_{1}%
,d_{2})=(0.4,2,3).$ Both graphs have the same average degree of $d=2.4$.
(c)\ Small-world network consisting of a ring and with edges added at random
between the odd-numbered vertices. It has an average degree of $d=2.5.$}%
\label{fig:graphs}%
\end{figure}

Numerous papers address various aspects of this question, see e.g.
\cite{olfati2005ultrafast,kolokolnikov2015maximizing, ogiwara2015maximizing,
shahbaz2021algebraic}. This is one of those situations where the best answer
is elusive \cite{mosk2008maximum}, but a decent answer can be found relatively
quickly \cite{ghosh2006growing, wang2008algebraic}. Graphs with high algebraic
connectivity are related to expander graphs, and are important in many
applications \cite{hoory2006expander, lubotzky2012expander}. In this paper, we
study the algebraic connectivity of \emph{sparse} random graphs in the case
when the number of edges $m$ scales linearly with the number of vertices $n,$
\ i.e. $m=O(n)$. In other words, we fix the average degree $d=2m/n=O(1)$ while
letting $n\rightarrow\infty.$

It is well known that a random Erdos-Renyi graph (where $m$ edges are taken at
random) needs $O(n\log n)$ edges to be fully connected \cite{erdos,
alon2004probabilistic}, and as such, they are poor candidates for maximizing
algebraic connectivity in this sparse regime (since AC\ of a disconnected
graph is zero). A good candidate are $d-$regular graphs for which every vertex
has degree $d$. It is well known that random regular $d-$graphs have an
algebraic connectivity that asymptotes to $\mu\sim d-2\sqrt{d-1}$ as
$n\rightarrow\infty$ with fixed $d$ \cite{mckay1981expected,
alon1986eigenvalues, broder1987second, friedman1991second, nilli2004tight}.
For integer $d\geq3,$ this is quantity is bounded away from zero, which
assures these graphs have good expander properties \cite{wormald1999models,
murty2003ramanujan}. A natural question is whether one can do better than
regular graph, but without \textquotedblleft too much work\textquotedblright%
\ (in, say, $O(n)$ time). In this paper we give an affirmative answer when
$d\leq7$:\ we introduce a class of semi-regular random graphs (whose vertices
have degree either $d_{1}$ or $d_{2}\neq d_{1}$) which are as easy to
construct as random $d$-regular graphs, but which have better AC.

Another question is whether \textquotedblleft expander-type\textquotedblright%
\ graphs (which we define to be graphs with AC bounded away from 0 as
$n\rightarrow\infty$)\ are possible when the average degree $d$ is less than
3, with $n\rightarrow\infty.$ For $2<d<3,$ the answer is yes, and it is
provided by semi-regular graphs of degrees $2$ and $3,$ whose average degree
is $2+p$ (with $0<p<1)$. We show that AC\ of such graphs asymptotes to
$\mu\sim p^{2}/4$ in the limit of small $p$ and large $n$, with $n\rightarrow
\infty$ independent of $p$. (The answer is no when $d\leq2$ since any graph
needs at least $n-1$ edges to be connected).

Our first result is on random semi-regular \emph{bipartite }graphs, where each
vertex has degree either $d_{1}$ or $d_{2}$, and we compute their asymptotic
AC in the limit $n\rightarrow\infty.$ In particular, we will exhibit a family
of semi-random bipartite graphs having the same average degree as random
$d$-regular graphs, but which have higher algebraic connectivity when
$d\leq7.$

Before stating our result, let us define what we mean by such graphs. Consider
a bipartite graph with one part having $n_{1}$ vertices of degree $d_{1},$ and
the second part having $n_{2}$ vertices of degree $d_{2},$ such that every
edge is between these two parts, and no edges are within each part. Then
$n_{1}d_{1}=n_{2}d_{2}.$ Such graph has $n=n_{1}+n_{2}$ vertices so that%
\begin{equation}
n_{1}=\frac{d_{2}}{d_{1}+d_{2}}n,\ \ \ n_{2}=\frac{d_{1}}{d_{1}+d_{2}}n
\end{equation}
and therefore the average degree is%
\begin{equation}
d=\frac{2d_{1}d_{2}}{d_{1}+d_{2}}. \label{d}%
\end{equation}
We call such a graph a $(d_{1},d_{2})$ semi-regular bipartite graph. We now
introduce the following random model.

\textbf{Random semi-regular bipartite (RSRB)\ graph model}:\ Take two bags. In
the first bag, put $d_{1}$ copies of $n_{1}$ vertices labelled $1\ldots
n_{1}.$ In the second bag, put $d_{2}$ copies of $n_{2}$ vertices labelled
$n_{1}+1\ldots n_{1}+n_{2}.$ Start with an empty graph of $n_{1}+n_{2}$
vertices. Then randomly pick two vertices without replacement -- one from each
bag -- and add an edge between them to the graph. Repeat until the bags are
empty. Refer to Figure \ref{fig:graphs}(a) and Matlab code in Appendix A.

We remark that this model (and the theory below) generally allows for multiple
edges. If desired, they can be eliminated through a random rewiring
postprocessing step\footnote{Given a multiple edge $(a,b)$ or a self loop
($a=b$), choose another edge $\left(  c,d\right)  $ at random. Then replace
edges $\left(  a,b\right)  $ and $\left(  c,d\right)  $ by edges $\left(
a,c\right)  $ and $\left(  b,d\right)  .$ Repeat until all multiple
edges/loops are eliminated. This operation preserves degree distribution and
edge count.}. These cases are sufficiently rare that the postprocessing step
does not effect the asymptotic results in the large $n$ limit. We now state
our main results for RSRB\ graphs.

\begin{mainres}
\label{thm:spectrum}Consider a $\left(  d_{1},d_{2}\right)  $ RSRB graph. In
the limit $n\rightarrow\infty,$ its spectrum density asymptotes to%
\begin{equation}
\rho(x)=\left\{
\begin{array}
[c]{c}%
\frac{1}{\pi}\frac{d_{1}d_{2}}{d_{1}+d_{2}}\frac{\sqrt{(x^{2}-r_{-}^{2}%
)(r_{+}^{2}-x^{2})}}{\left(  d_{1}d_{2}-x^{2}\right)  \left\vert x\right\vert
},\ \ \left\vert x\right\vert \in\left(  r_{-},r_{+}\right) \\
\frac{\left\vert d_{2}-d_{1}\right\vert }{d_{1}+d_{2}}\delta\left(  x\right)
,\ \ \left\vert x\right\vert <r_{-}\\
0,\ \ \ \left\vert x\right\vert >r_{+}%
\end{array}
\right.  \label{rho}%
\end{equation}
where $\delta$ is the Dirac-delta function and
\begin{equation}
r_{\pm}=\left(  d_{1}+d_{2}-2\pm\sqrt{\left(  d_{1}+d_{2}-2\right)
^{2}-\left(  d_{2}-d_{1}\right)  ^{2}}\right)  ^{1/2}. \label{r}%
\end{equation}
In other words, the number of eigenvalues inside any interval $\left(
a,b\right)  $ asymptotes to $\int_{a}^{b}\rho(x)dx$ as $n\rightarrow\infty.$

Moreover, its algebraic connectivity asymptotes to%
\begin{equation}
\mu\sim\frac{d_{1}+d_{2}}{2}-\left(  \left(  \frac{d_{2}-d_{1}}{2}\right)
^{2}+r_{+}^{2}\right)  ^{1/2},\ \ \ n\gg1 \label{muRSRB}%
\end{equation}

\end{mainres}

\begin{figure}[tbh]
\includegraphics[width=1\textwidth]{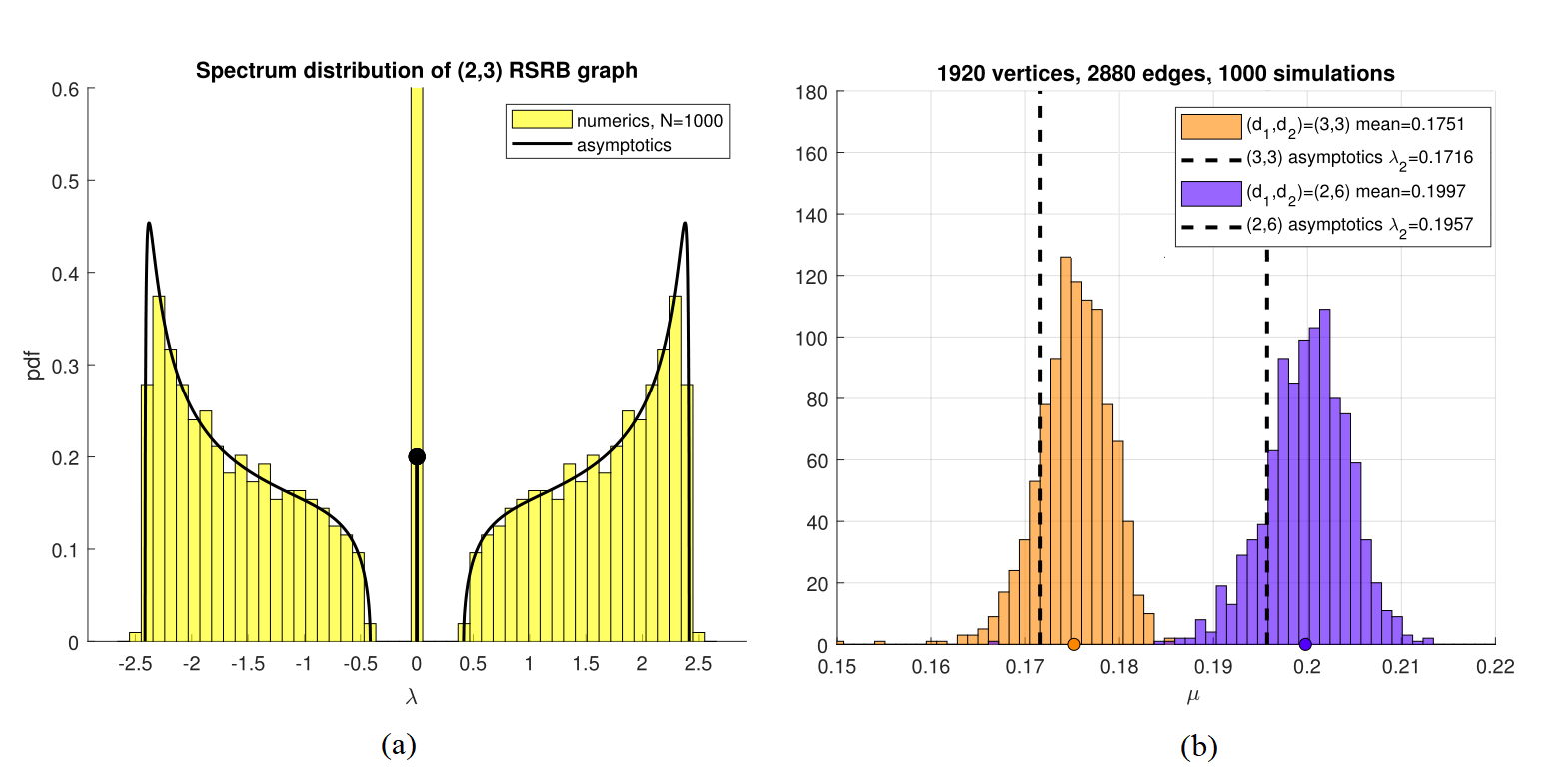}\caption{ (a) Full spectrum
of random semi-regular bipartite graph with $\left(  d_{1},d_{2}\right)
=\left(  2,3\right)  $. Numerics correspond to the histograph of eigenvalues
of a single such graph with 1000 vertices, computed numerically using Matlab.
Asymptotics corresponds to the formula (\ref{rho}). The height of the lollypop
corresponds to the weight delta function at the origin. (b) Comparison of
algebraic connectivity between $\left(  3,3\right)  $ regular bipartite,
$\left(  2,6\right)  $ semi-regiular bipartite graphs, and the asymptotic
theory. The two classes have the same number of vertices and edges, and
$\left(  2,6\right)  $ is 15\% better than (3,3) (both for asymptotics and
numerics).}%
\label{fig:rsrb}%
\end{figure}

Figure \ref{fig:rsrb}(a) shows the shape of the distribution (\ref{rho})\ for
the case $d_{1}=2,\ \ d_{2}=3$ (having an average degree of $d=2.4$), and
compares it to an numerical histogram of eigenvalues of a 1000-vertex graph
for this case. Very good agreement is observed.

Note that in the case $d_{1}=d_{2},$ formula (\ref{rho})\ reduces to the
well-known Mackay distribution of the spectrum of a regular graph
\cite{mckay1981expected}:
\begin{equation}
\rho(x)=\left\{
\begin{array}
[c]{c}%
\frac{1}{\pi}\frac{d}{2}\frac{\sqrt{4(d-1)-x^{2}}}{d^{2}-x^{2}},\ \ \left\vert
x\right\vert <2\sqrt{d-1}\\
0,\ \ \left\vert x\right\vert >2\sqrt{d-1}.
\end{array}
\right.
\end{equation}
Moreover, formula (\ref{muRSRB})\ simplifies to the classical result $\mu\sim
d-2\sqrt{d-1}$ for the algebraic connectivity of a random regular graph of $d$
vertices \footnote{This also shows that a random regular graph and random
bipartite regular graphs have the same algebraic connectivity, at least to
leading order in $n.$}.

Consider a random cubic graph $\left(  d_{1}=d_{2}=3\right)  $ and contrast it
with an RSRB\ graph with $d_{1}=2,\ d_{2}=6.$ Both cases have the same average
degree of $d=3.$ Formula (\ref{muRSRB})\ gives the asymptotic values of
$\mu\sim0.17157$ for the former, and $\mu\sim0.19577,$ for the latter.
\emph{Thus, }$\left(  2,6\right)  $\emph{ RSRB\ graph is about 15\% better
than a random cubic graph with respect to its expander properties, while
having the same number of vertices and edges.} Figure \ref{fig:rsrb}(b) shows
a histogram for $\mu$ of 1000 randomly constructed such graphs, comparing
these two cases with $n=2000$. Very good agreement between numerics and
asymptotics is observed in both cases.

More generally, the following table shows all possible combinations of
$d_{1},d_{2}$ such that $d=\frac{2d_{1}d_{2}}{d_{1}+d_{2}}$ is an integer
between 3 and 8, and the corresponding value of $\mu$ (\ref{muRSRB}).

\begin{center}%
\begin{tabular}
[c]{|c|c|c|c|c|c|c|c|c|c|c|c|c|c|}\hline
\multicolumn{14}{|c|}{All RSRB graphs with integer average degree
$d=3,\ldots,8$}\\\hline
$d$ & \multicolumn{2}{|c|}{3} & \multicolumn{2}{|c|}{4} &
\multicolumn{2}{|c|}{5} & \multicolumn{2}{|c|}{6} & \multicolumn{2}{|c|}{7} &
\multicolumn{3}{|c|}{8}\\\hline
$d_{1}$ & 3 & 2 & 4 & 3 & 5 & 3 & 6 & 4 & 7 & 4 & 8 & 6 & 5\\\hline
$d_{2}$ & 3 & 6 & 4 & 6 & \ 5 \  & 15 & 6 & 12 & 7 & 28 & 8 & 12 & 20\\\hline
$\mu_{\text{asympt}}$ & 0.1715 & \textbf{0.1957} & 0.5358 & \textbf{0.5535} &
1 & \textbf{1.0890} & 1.5278 & \textbf{1.5587} & 2.1010 & \textbf{2.1435} &
\textbf{2.7084} & 2.6887 & 2.6671\\\hline\hline
$\mu_{\text{numerics}}$ & 0.178 & 0.205 & 0.553 & 0.572 & 1.027 & 1.122 &
1.565 & 1.596 & 2.150 & 2.205 & 2.766 & 2.745 & 2.729\\\hline
std & 0.006 & 0.006 & 0.011 & 0.010 & 0.015 & 0.017 & 0.018 & 0.018 & 0.021 &
0.020 & 0.026 & 0.022 & 0.022\\\hline
diff \% & 3.8\% & 4.7\% & 3.1\% & 3.2\% & 2.7\% & 3.0\% & 2.4\% & 2.4\% &
2.3\% & 2.87\% & 2.1\% & 2.1\% & 2.2\%\\\hline
\end{tabular}

\end{center}

The row $\mu_{\text{asympt}}$ is the asymptotic formula given
by\ (\ref{muRSRB}). Row $\mu_{\text{numerics}}$\ corresponds to Monte Carlo
simulations of $\mu$. It shows the average $\mu$ for 200 randomly chosen
RSRB\ graphs with $n=1000$ edges. The row \textquotedblleft diff
\%\textquotedblright\ is $\frac{\mu_{\text{numerics}}-\mu_{\text{asympt}}}%
{\mu_{\text{asympt}}}\times100$. Uniformly good agreement with between
asymptotics and numerics is observed.

Parameters with higher AC\ are in shown in bold. For $d\leq7,$ RSRB with
$d_{1}\neq d_{2}$ have higher algebraic connectivity than the $d$-regular
graph. On the other hand, for $d\geq8$, $d$-regular graphs win.

It is interesting to note that for any integer $d\geq3,$ equation (\ref{d})
always has a solution with integers $2\leq d_{1}<d_{2}.$ When $d$ is prime,
this solution is unique and is given by $d_{1}=(d+1)/2,$ $d_{2}=d\left(
d+1\right)  /2.$\ More generally, the number of such solutions is precisely
the number of pythegorean triples of leg $d$ (sequence A046079 in OEIS).

RSRB graphs above have a constraint $n_{1}d_{1}=n_{2}d_{2}.$ In particular the
minimum attainable average degree of such graphs is $d=2.4$ corresponding to
$\left(  d_{1},d_{2}\right)  =\left(  2,3\right)  .$ We can remove this
constraint by instead introducing the probability of having degree $d_{1}$ or
$d_{2}$ as follows.

\textbf{Random semi-regular (RSR)\ model:} Given $p,d_{1},d_{2},$ and $n,$ let
$n_{1}=\left\lfloor \left(  1-p\right)  n\right\rfloor $ and let
$n_{2}=n-n_{1}.$ In the same bag, put $d_{1}$ copies of $n_{1}$ vertices
labelled $1\ldots n_{1}$ and $d_{2}$ copies of $n_{2}$ vertices labelled
$n_{1}+1\ldots n.$ Create edges by drawing two vertices from the bag at random
(without replacement), until only one or zero vertices are left in the bag.
See Appendix A for Matlab code.

An example of RSR\ graph is shown in Figure \ref{fig:graphs}(b). Note that
such a graph has an average degree of $d=\left(  1-p\right)  d_{1}+pd_{2}$. We
have the following.

\begin{mainres}
\label{main:2}Consider a $\left(  p,\ d_{1},d_{2}\right)  $ random
semi-regular graph. Let
\begin{equation}
F(R,x)=x\left(  d_{{2}}-d_{{1}}\right)  \left(  1-Rx\right)  p+\left(
R{x}^{2}(d_{{2}}-1)-1\right)  \left(  {R}^{2}{x}^{2}\left(  d_{{1}}-1\right)
+Rx\left(  d_{{2}}-d_{1}\right)  -R+1\right)  . \label{F(R,x)}%
\end{equation}

Let $x$ be the smallest root of the system $F=0=\partial F/\partial R.$Then in
the limit $n\rightarrow\infty,$ the AC is given by $\mu=d_{2}-1/x.$
\end{mainres}

In general, eliminating $R$ from the system $F=0=\partial F/\partial R$ is a
straightforward computer algebra computation using a resultant, and yields in
a 6th degree polynomial for $x$. It is too ugly to write down here for general
$d_{1},d_{2}$ -- see Appendix A for Maple code. In the case $d_{1}=2,d_{2}=3,$
RSR graph has average degree $2+p,$ and $\mu$ is the smallest root of%
\begin{equation}
0=\mu\left(  \mu-4\right)  \left(  \mu^{2}-4\mu-1\right)  +2\mu\left(
3\mu^{3}-33\mu^{2}+89\mu-19\right)  p+\left(  -15\mu^{2}-30\mu+1\right)
p^{2}+8p^{3}. \label{mu23}%
\end{equation}

\begin{figure}[tb]
\includegraphics[width=0.6\textwidth]{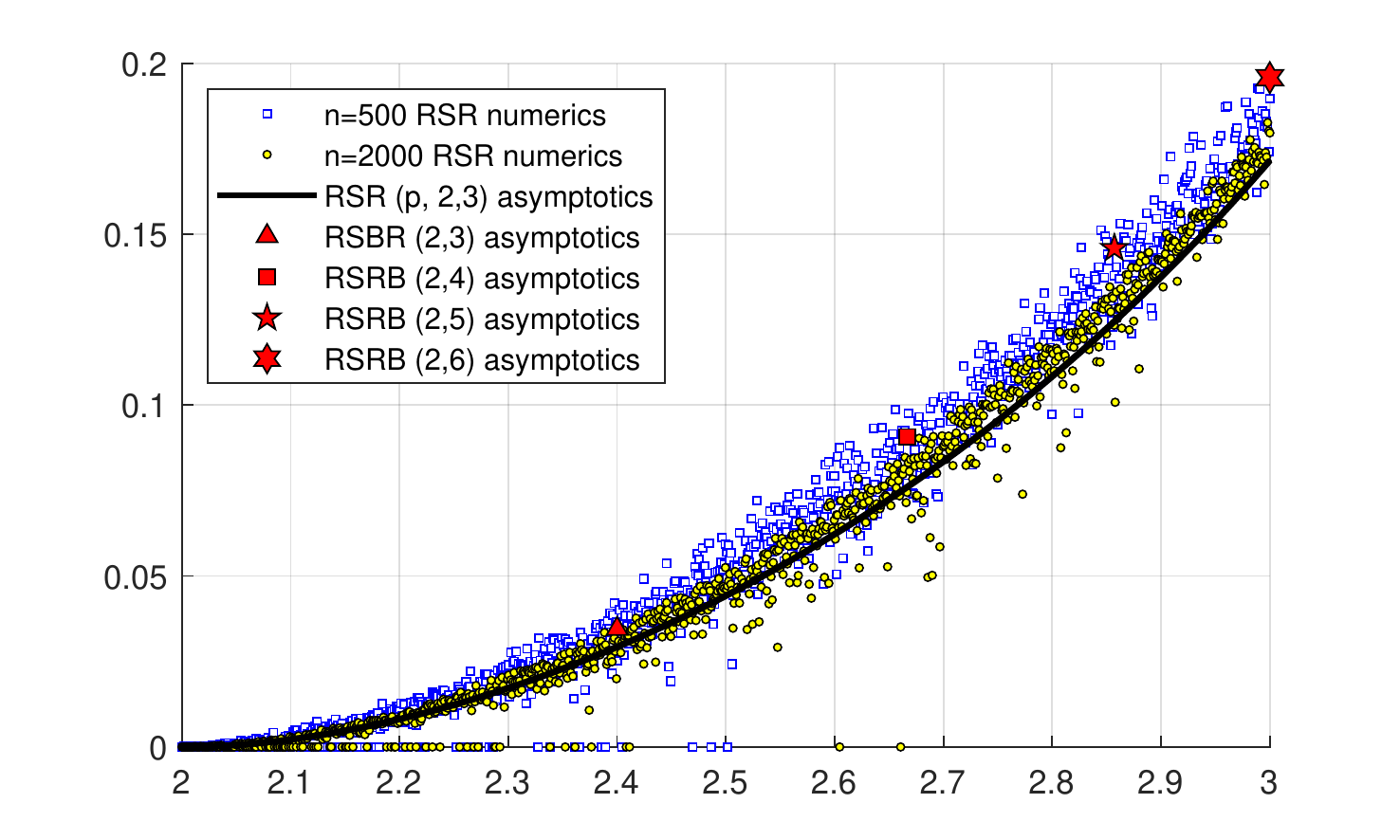}\caption{Comparison between
RSR, RSRB asymptotics, and RSR numerics, with average degree $2<d\leq3.$
Numerics represent AC of 1000 randomly chosen RSR graphs with $p\in(0,1)$ and
$(d_{1},d_{2})=(2,3)$. Here, $\mu$ is plotted against $d=2+p$. Asymptotics for
RSR correspond to roots of (\ref{mu23}). Asymptotics for RSRB are given by
(\ref{muRSRB}).}%
\label{fig:RSR}%
\end{figure}

Figure \ref{fig:RSR} compares $\mu$ given by (\ref{mu23}) with numerical
computations of $\mu$ for randomly chosen $\left(  p,2,3\right)  $
RSR\ graphs. Note that the numerical result approaches the asymptotic value of
$\mu$ as the number of edges $n$ is increased.\bigskip

The following table gives the value of $\mu$ for several choices of
$p,d_{1},d_{2}$ for which the average degree $d=4,$ and comparison with numerics:%

\[%
\begin{tabular}
[c]{|c|c|c|c|c|c|c|}\hline
\multicolumn{7}{|c|}{$\text{RSR\ model with }d=4$}\\\hline
$d_{1}$ & 4 & 3 & 3 & 2 & 2 & 2\\\hline
$d_{2}$ & 4 & 5 & 6 & 5 & 6 & 7\\\hline
$p$ &  & 0.5 & 1/3 & 2/3 & 0.5 & 0.4\\\hline
\multicolumn{1}{|c|}{$\mu_{\text{asympt}}$} & \textbf{0.5359} & 0.44261 &
0.39162 & 0.3333 & 0.25352 & 0.20748\\\hline\hline
\multicolumn{1}{|c|}{$\mu_{\text{numerics}}$} & 0.551 & 0.488 & 0.451 &
0.286 & 0.217 & 0.174\\\hline
std & 0.010 & 0.020 & 0.022 & 0.062 & 0.051 & 0.045\\\hline
diff \% & 2.8\% & 10\% & 15\% & -14\% & -14\% & -19\%\\\hline
\end{tabular}
\ \
\]
Here, $\mu_{\text{asympt}}$ is as given by Main Result \ref{main:2};
$\mu_{\text{numerics}}$ is the average of 200 simulations with $n=1000.$ As
apparent from this table, RSR model performs worse than RSRB of average degree
4 (including $4-$regular graphs). This appears to be true for any average
degree $d$. However one advantage of RSR\ graphs is that they produce graphs
with an average degree $2<d<3,$ with $\mu$ bounded away from zero as
$n\rightarrow\infty.$ The smallest average degree that RSRB model can have is
2.4 corresponding to $\left(  d_{1},d_{2}\right)  =\left(  2,3\right)  .$ In
fact, (\ref{mu23}) shows that $\mu\sim p^{2}/4$ when the average degree is
$2+p$ with $d_{1}=2,d_{2}=3,$ and $0<p\ll1.$

\section{Derivation of main results}

\textbf{Derivation of spectral density \ (\ref{rho}).} Following \cite{wigner,
mckay1981expected} we use the trace method. It is based on the fact that
$\operatorname{trace}(A^{s})=\sum_{j=1}^{s}\lambda_{j}^{s}.$Define $\phi
_{s}=\frac{1}{n}\operatorname{trace}(A^{s}).$ Then $\phi_{s}$ represents the
average number of closed walks of length $s$ on the graph whose adjacency
matrix is $A$. In the limit $n\rightarrow\infty,$ the eigenvalue distribution
therefore satisfies%
\begin{equation}
\int x^{s}\rho(x)dx=\phi_{s}. \label{1206}%
\end{equation}
As in \cite{wigner, mckay1981expected}, the derivation consists of (a)
computing $\phi_{s};$ and (b)\ inverting the integral equation (\ref{1206}%
)\ to determine $\rho(x).$

To compute $\phi_{s},$ we write%
\[
\phi_{s}=\frac{d_{1}}{d_{1}+d_{2}}\phi_{A,s}+\frac{d_{2}}{d_{1}+d_{2}}%
\phi_{B,s},
\]
where $\phi_{A,s}$ is the number of closed walks of length $s$ starting from a
vertex of degree $d_{2},$ and $\phi_{B,s}$ is that number for vertices of
degree $d_{1}.$ A key insight \cite{mckay1981expected} that allows the
asymptotic determination of $\phi_{A,s}$ and $\phi_{B,s}$ is that locally, a
large random regular graph looks like a tree because the probability of
encountering a cycle of length $s$ is vanishinly small as $n\rightarrow\infty$
(for fixed $s$). The same is true of semi-regular random graphs. For the RSRB
graphs, each successive level alternates between vertices of degree $d_{1}$
and $d_{2}.$ We therefore decompose these trees as illustrated in the
following diagram, for the case $d_{1}=2$ and $d_{2}=3$:

\begin{center}%
\begin{gather*}
d_{1}=2,\ \ d_{2}=3\\
\includegraphics[width=1\textwidth]{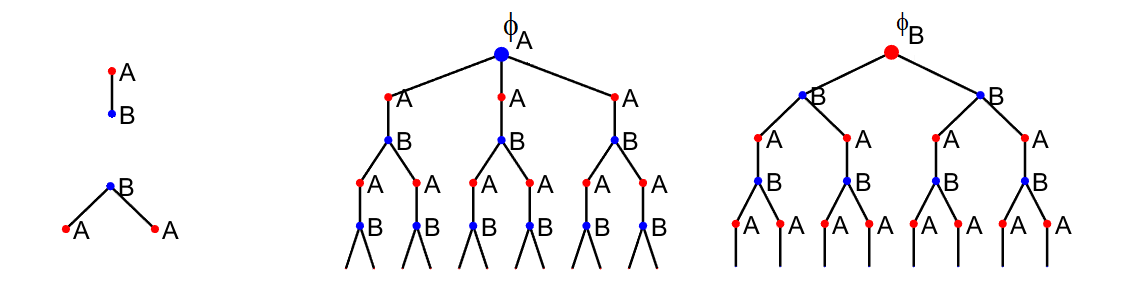}
\end{gather*}

\end{center}

Accordingly, $\phi_{A,s}$ and $\phi_{B,s}$ satisfy the recursion%
\begin{align}
\phi_{A,s}  &  =d_{2}\sum_{j=0}^{s-2}\phi_{A,j}A_{s-2-j},\ \ \ \ \phi
_{B,s}=d_{1}\sum_{j=0}^{s-2}\phi_{B,j}B_{s-2-j},\nonumber\\
A_{s}  &  =\left(  d_{1}-1\right)  \sum_{j=0}^{s-2}A_{j}B_{s-2-j}%
,\ \ \ \ \ B_{s}=\left(  d_{2}-1\right)  \sum_{j=0}^{s-2}B_{j}A_{s-2-j},
\label{1208}%
\end{align}
with $A_{0}=B_{0}=\phi_{A,0}=\phi_{B,0}=1.$ Let $A(x)=\sum A_{s}x^{s}$ be the
generating function for $A_{s}$ and similarly for $B,\ \phi_{A},\ \phi
_{B},\ \phi.$ From recursion relations (\ref{1208}), the corresponding
generating functions satisfy%
\begin{align}
A  &  =1+\left(  d_{1}-1\right)  x^{2}AB,\nonumber\\
B  &  =1+\left(  d_{2}-1\right)  x^{2}BA,\nonumber\\
\phi_{A}  &  =1+d_{2}x^{2}\phi_{A}A,\label{1213}\\
\phi_{B}  &  =1+d_{1}x^{2}\phi_{B}B,\nonumber\\
\phi &  =\frac{d_{1}}{d_{1}+d_{2}}\phi_{A}+\frac{d_{2}}{d_{1}+d_{2}}\phi
_{B}\nonumber
\end{align}

Solving (\ref{1213}), while keeping in mind that $\phi\left(  0\right)  =1,$we
find that%
\[
\phi(x)=\frac{d_{1}d_{2}}{d_{1}+d_{2}}\frac{\frac{d_{1}+d_{2}}{d_{1}d_{2}%
}-\sqrt{\left(  d_{2}-d_{1}\right)  ^{2}x^{4}+\left(  4-2d_{1}-2d_{2}\right)
x^{2}+1}-1}{d_{1}d_{2}x^{2}-1}.
\]
In principle, $\phi_{s}$ can be computed from $\phi(x)$ by Taylor-expanding
near the origin; the first few terms are $\phi_{0}=1,\ \phi_{2}=2\frac
{d_{1}d_{2}}{d_{1}+d_{2}},\ \phi_{4}=2\frac{d_{1}d_{2}}{d_{1}+d_{2}}\left(
d_{1}+d_{2}-1\right)  $. In \cite{mckay1981expected} explicit computation of
$\phi_{s}$ was combined with Tchebychev polynomials to compute the McKay
distribution for $d$-regular random graphs. Here, we use a simpler and more
powerful technique using complex variables introduced in
\cite{longoria2020pascal}.

Write $\phi_{s}$ using Cauchy's integral formula as
\[
\phi_{s}=\int z^{-s-1}\phi(z)\frac{dz}{2\pi i},
\]
where the integration is around a circle $\left\vert z\right\vert
=\varepsilon$, with $\varepsilon$ sufficiently small to avoid any branch cuts
of $\phi.$ Using the fact that $\phi(z)$ is even, we then have%

\begin{align*}
\phi_{s}  &  =\frac{1}{\pi}\int_{0}^{\pi}\varepsilon^{-s}\phi\left(
\varepsilon e^{i\theta}\right)  e^{-si\theta}d\theta\\
&  =\frac{i}{\pi}\int_{-1/\varepsilon}^{1/\varepsilon}\frac{1}{x}\phi\left(
\frac{1}{x}\right)  x^{s}dx
\end{align*}
Taking $\varepsilon\rightarrow0$ and recalling (\ref{1206})\ we obtain%
\[
\int\rho(x)x^{s}dx=\frac{i}{\pi}\int_{-\infty}^{\infty}\frac{1}{x}\phi\left(
\frac{1}{x}\right)  x^{s}dx.
\]
This yields the formula for $\rho,$ namely%
\begin{equation}
\rho(x)=-\frac{1}{\pi}\operatorname{Im}\left\{  \frac{1}{x}\phi\left(
\frac{1}{x}\right)  \right\}  .
\end{equation}
This formula was derived in \cite{longoria2020pascal} using an equivalent
technique (Stieltjes inversion). Next we compute%
\begin{equation}
\frac{1}{x}\phi\left(  \frac{1}{x}\right)  =\frac{1}{x}\frac{d_{1}d_{2}}%
{d_{1}+d_{2}}\frac{\frac{d_{1}+d_{2}}{d_{1}d_{2}}x^{2}-\sqrt{(x^{2}-r_{-}%
^{2})(x^{2}-r_{+}^{2})}-x^{2}}{d_{1}d_{2}-x^{2}}, \label{G}%
\end{equation}
where $r_{\pm}$ are given by\ (\ref{r}). It follows that $\operatorname{Im}%
\left\{  \frac{1}{x}\phi\left(  \frac{1}{x}\right)  \right\}  =0$ for
$\left\vert x\right\vert \notin\left(  r_{-},r_{+}\right)  $ and $x\neq0.$ On
the other hand when $\left\vert x\right\vert \in\left(  r_{-},r_{+}\right)  $
we obtain
\begin{equation}
\rho(x)=\frac{1}{\pi}\frac{1}{\left\vert x\right\vert }\frac{d_{1}d_{2}}%
{d_{1}+d_{2}}\frac{\sqrt{(r_{+}^{2}-x^{2})(x^{2}-r_{-}^{2})}}{d_{1}d_{2}%
-x^{2}},\ \ \left\vert x\right\vert \in\left(  r_{-},r_{+}\right)  .
\label{1524}%
\end{equation}
This recovers the formula (\ref{rho})\ except for the delta mass at the center
which is due to the singularity there. To compute the weight of the delta
mass, one can integrate the overall density and impose that $\int\rho=1.$
While the exact integration of (\ref{1524}) is possible, it is easier and more
instructive to compute the number of zero eigenvalues directly as follows.

Suppose that $d_{2}>d_{1}.$ Label the components of the eigenvector
$x_{1}\ldots x_{n_{1}}$ and $y_{1}\ldots y_{n_{2}},$ corresponding to vertices
of degree $d_{1}$ and $d_{2},$ respectively. A zero eigenvalue satisfies
$\sum_{x_{j}\in nbr(y_{k})}x_{j}=0$,\ \ $k=1\ldots n_{2}$ and $\sum_{y_{j}\in
nbr(x_{k})}y_{j}=0$, \ \ $k=1\ldots n_{1}.$ Look for solutions of the form
where $y_{j}=0$, $j=1\ldots n_{2}.$ This corresponds to solving $n_{2}$ linear
equations $\sum_{x_{j}\in nbr(y_{k})}x_{j}=0$, for the $n_{1}$ unknows
$x_{1}\ldots x_{n_{1}}.$ This is an under-determined system, since
$n_{2}<n_{1}.$ Generically, it has $n_{1}-n_{2}$ independent solutions.
Therefore the weight of the delta function at zero should be $\frac
{n_{1}-n_{2}}{n}.$ Using (\ref{d})\ yields $\frac{n_{1}-n_{2}}{n}=\frac
{d_{2}-d_{1}}{d_{1}+d_{2}},$ which recovers the weight of the delta function
in (\ref{rho}). We also verified using computer algebra that this indeed
agrees with the full integration of (\ref{1524}), namely, that $\int_{r_{-}%
}^{r_{+}}\rho(x)dx=\frac{d_{1}}{d_{1}+d_{2}}$, \ \ \ $d_{2}>d_{1},$ so that
indeed $\int_{-\infty}^{\infty}\rho=1$. This completes the derivation of the
spectral density (\ref{rho}).

Finally, another derivation of the delta weight, as pointed out by James
Mingo, is to use Proposition 3.8 in \cite{mingo2017free}. It says that if
$\frac{1}{x}\phi\left(  \frac{1}{x}\right)  $ has a singularity at $x=a,$ then
the associated measure $\rho(x)$ has a delta mass at $x=a$ whose weight is
given by $\lim_{x\rightarrow a}\left\vert \left(  x-a\right)  \frac{1}{x}%
\phi\left(  \frac{1}{x}\right)  \right\vert .$ Here, $a=0$ and the mass weight
is then given by $\lim_{x\rightarrow0}\left\vert \phi\left(  \frac{1}%
{x}\right)  \right\vert =\frac{d_{1}d_{2}}{d_{1}+d_{2}}\frac{\sqrt{r_{-}%
^{2}r_{+}^{2}}}{d_{1}d_{2}}=\frac{\left\vert d_{2}-d_{1}\right\vert }%
{d_{1}+d_{2}}$ in agreement with direct computation. $\blacksquare$

\textbf{Derivation of AC\ formula for RSRB\ }(\ref{muRSRB}) Next we derive the
formula for algebraic connectivity (\ref{muRSRB}). For $d$-regular graphs, the
Laplacian graph is given by $dI-A$ where $I$ is the identity and $A$ is the
adjacency matrix. This allows to characterize the AC in terms of
second-largest eigenvalue $\lambda$ of $A$:\ $\mu=d-\lambda.$ Of course this
is not true if the graph is not regular. The trick is to \textbf{regularize
the graph by adding enough self-loops to vertices until all vertices have the
same degree}. Each self-loops counts as a single additional edge (i.e. adding
a loop to vertex $j$ corresponds to adding "1" to the $j$-th diagonal entry of
the associated adjacency matrix) and crucially, adding self-loops does not
change the Laplacian of the graph.

For a $\left(  d_{1},d_{2}\right)  $ semi-regular graph, assume that
$d_{1}<d_{2}$ and add $d_{2}-d_{1}$ loops to all vertices of degree $d_{1}.$
This results in a $d_{2}$-regular graph with loops. As before let $\phi_{s}$
be the number of closed walks of size $s$. Let $\lambda$ be the \emph{local
expansion rate}, that is, the rate at which $\phi_{s}$ grows: $\lambda
=\lim\limits_{s\rightarrow\infty}\phi_{s+1}/\phi_{s}$. Then by analogy to
regular graphs, the AC is given by $\mu=d_{2}-\lambda.\,\ $To compute
$\lambda,$ we decompose closed loops similarly to (\ref{1213}). The
decomposition now has loops as illustrated below:%

\begin{gather*}
d_{1}=2,\ d_{2}=3\\
\includegraphics[width=0.95\textwidth]{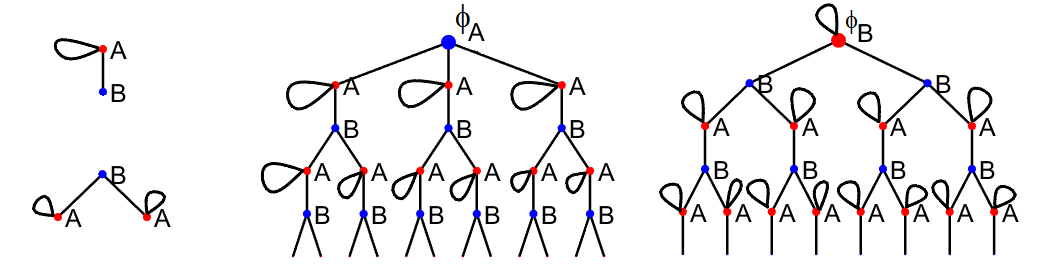}
\end{gather*}

Correspondingly, the associated generating functions satisfy%
\begin{align}
A  &  =1+(d_{2}-d_{1})xA+\left(  d_{1}-1\right)  x^{2}AB,\ \ \ \nonumber\\
B  &  =1+\left(  d_{2}-1\right)  x^{2}BA,\nonumber\\
\phi_{A}  &  =1+d_{2}x^{2}\phi_{A}A,\label{119}\\
\phi_{B}  &  =1+x\left(  d_{2}-d_{1}\right)  \phi_{B}+d_{1}x^{2}\phi
_{B}B,\nonumber\\
\phi &  =\frac{d_{1}}{d_{1}+d_{2}}\phi_{A}+\frac{d_{2}}{d_{1}+d_{2}}\phi
_{B}\nonumber
\end{align}

To determine AC, it suffices to determine the growth rate of $\phi_{s}$. This
growth rate is in turn determined by the locations of singularities in the
associated generating functions \cite{wilf2005generatingfunctionology}. In
particular, if $\sum\phi_{s}z^{s}=g(z)+h(z)(z-r)^{p}$, where $p$ is
non-integer and $g,h$ analytic, then $\phi_{s}$ grows with the rate
$r^{-1}:\ $in other wors, $\phi_{s+1}/\phi_{s}\rightarrow r^{-1}$ as
$s\rightarrow\infty.$ The location of the singularity can be found without
solving the full system, and corresponds to the smallest root of the
discriminant of the system. To illustrate this, consider the Catalan numbers
$c_{n}=\frac{1}{n+1}\left(
\begin{array}
[c]{c}%
2n\\
n
\end{array}
\right)  $ whose generating function satisfies $c=1+xc^{2}.$ Explicitly,
$c(x)=\frac{1-\sqrt{1-4x}}{2x},$ and has a singularity at $x=1/4,$
corresponding to the zero of the discriminant of the quadratic $xc^{2}-c+1=0.$
It follows that $c_{n+1}/c_{n}\rightarrow4$ as $n\rightarrow\infty$ (indeed
more precise asymptotics $c_{n}\sim\frac{1}{\sqrt{\pi}}4^{n}n^{-3/2}$ can also
be derived from its generating function, although here, we only need the
growth rate).

Let us now get back to the system (\ref{119}). Eliminating $B,$ we find that
$A$ satisfies a quadratic,
\[
\left(  d_{2}-1\right)  \left(  x\left(  d_{2}-d_{1}\right)  -1\right)
x^{2}A^{2}+\left[  x^{2}\left(  d_{2}-d_{1}\right)  -x\left(  d_{2}%
-d_{1}\right)  +1\right]  A-1=0.
\]
Setting its discriminant to zero, we obtain a quartic polynomial for $x$,
namely%
\begin{equation}
\left(  {d}_{1}-{d}_{2}\right)  ^{2}{x}^{4}+2\left(  {d}_{1}+{d}_{2}-2\right)
\left(  {d}_{2}-d_{1}\right)  {x}^{3}+\left(  (d_{1}-d_{2})^{2}-2\left(
{d}_{1}+{d}_{2}-2\right)  \right)  {x}^{2}+2\left(  {d}_{1}-{d}_{2}\right)
x-1=0. \label{120}%
\end{equation}

The expansion rate is then $\lambda=1/x,$ so that $\mu=d_{2}-1/x$.
Substituting $x=1/(\mu-d_{2})$ into (\ref{120}), one obtains a 4th degree
polynomial for $\mu,$%
\[
0=\mu^{4}-2\left(  d_{1}+d_{2}\right)  \mu^{3}+\left(  \left(  d_{1}%
+d_{2}\right)  ^{2}+2d_{1}d_{2}-2d_{1}-2d_{2}+4\right)  \mu^{2}+2\left(
d_{1}+d_{2}\right)  \left(  d_{1}+d_{2}-2-d_{1}d_{2}\right)  \mu+\left(
d_{1}d_{2}-d_{1}-d_{2}\right)  ^{2}%
\]

Shift this polynomial to eliminate the $\mu^{3}$ term by substituting
$\mu=y+\frac{d_{1}+d_{2}}{2}.$ Magically, it simultaneously elimintes the
$\mu^{1}$ term resulting in a quadratic for $y^{2}$:%

\[
0=y^{4}+\left(  4-2\left(  d_{1}+d_{2}\right)  -\left(  d_{2}-d_{1}\right)
^{2}/2\right)  y^{2}+\frac{\left(  d_{2}-d_{1}\right)  ^{2}\left(  d_{1}%
+d_{2}\right)  }{2}+\frac{\left(  d_{2}-d_{1}\right)  ^{4}}{16}%
\]
Solving for $y$ using quadratic formula yields the solution (\ref{muRSRB}).
$\blacksquare$

\textbf{Derivation of Result \ref{main:2}. }Again, we compute the average
number of closed walks of length $s$, $\phi_{s}.$ As before, each vertex looks
like a tree locally. Each child in this tree has probability $\left(
1-p\right)  $ chance of having degree $d_{1},$ and probability $p$ of having
degree $d_{2}.$ As in the derivation of (\ref{muRSRB}), assume that $d_{1}\leq
d_{2}$, and add $d_{2}-d_{1}$ self-loops to each vertex of degree $d_{1}.$ The
walks on these trees then decompose according to the following diagram (here
given in the case $d_{1}=2,\ d_{2}=3$):%

\begin{gather*}
d_{1}=2,\ \ d_{2}=3\\
\includegraphics[width=0.51\textwidth]{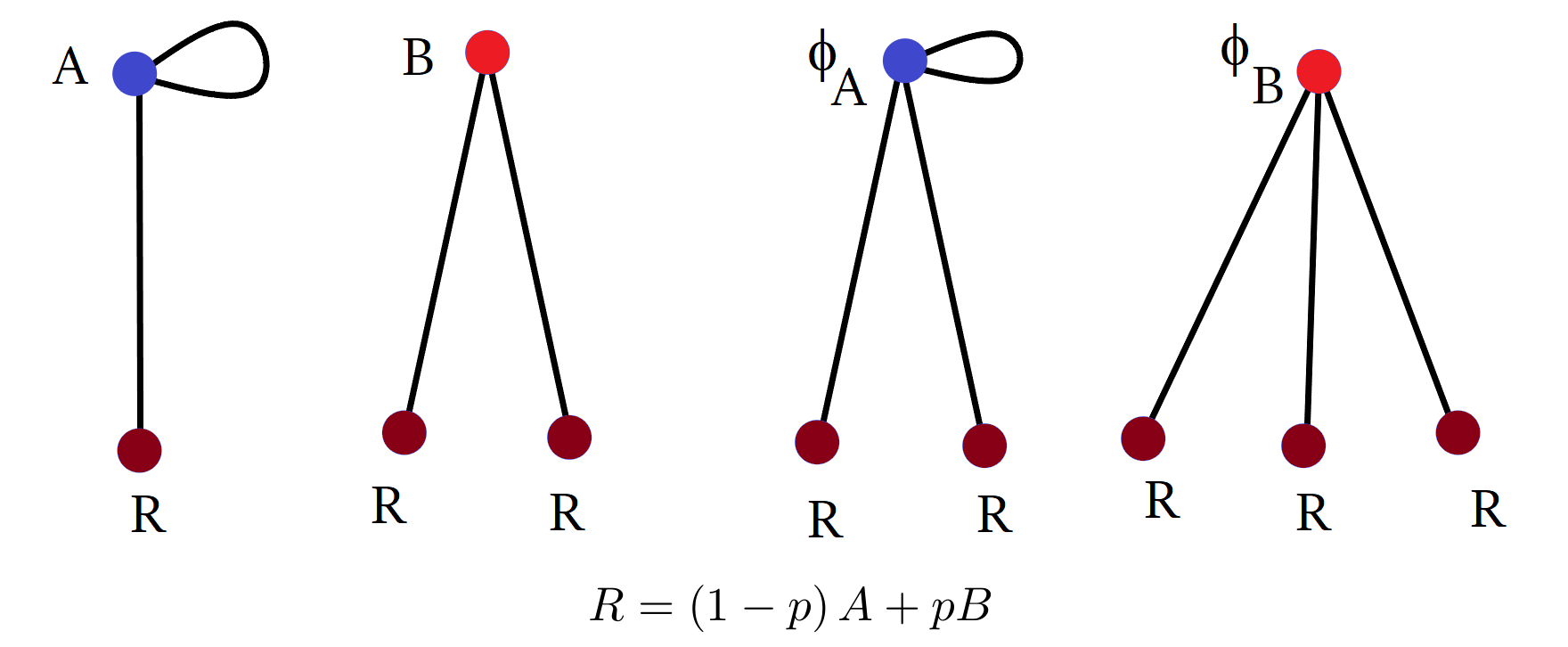}
\end{gather*}

The corresponding generating function $\phi(x)=\sum\phi_{s}x^{s}$ solves the equations%

\begin{align}
A  &  =1+x\left(  d_{2}-d_{1}\right)  A+\left(  d_{1}-1\right)  x^{2}%
AR\nonumber\\
B  &  =1+\left(  d_{2}-1\right)  BR\nonumber\\
R  &  =\left(  1-p\right)  A+pB\nonumber\\
\phi_{A}  &  =1+x\left(  d_{2}-d_{1}\right)  \phi_{A}+d_{1}x^{2}\phi_{A}R\\
\phi_{B}  &  =1+d_{2}x^{2}\phi_{B}R\nonumber\\
\phi &  =(1-p)\phi_{A}+p\phi_{B}.\nonumber
\end{align}

Eliminating $A$ and $B$ yields a cubic $F(R;x)=0$ given by (\ref{F(R,x)}). The
AC\ is then given by $\mu=d_{2}-1/x$, where $x$ is the singularity of $R(x)$
that is closest to the origin. By implicit function theorem, this happens when
$F_{R}=0.$ In other words, $x$ satisfies $F=0=\partial F/\partial
R.\blacksquare$

\section{Discussion and open problems}

We computed asymptotics of AC\ for two models of semiregular random
graphs:\ either RSR\ or RSRB\ models. While RSR\ model is shown to have
smaller AC than a random $d$-regular graph with the same average degree, we
found that RSRB model has \emph{higher}\textbf{ }AC than $d$-regular when
$d\leq7.$

Ramunajan expander graphs are defined \cite{hoory2006expander,
lubotzky2012expander} as being $d$-regular graphs whose algebraic connectivity
is at least $d-2\sqrt{d-1}$. As illustrated in Figure \ref{fig:rsrb}(b)\ (see
also Figure 9 in \cite{hoory2006expander} and related table there), a random
$d-$regular graph has a decent chance of being Ramunajan graph (around 66\%
when $d=3$ as $n\rightarrow\infty$ according to simulations in
\cite{hoory2006expander}). We can extend the definition of Ramunajan exander
graphs as being any graph of \emph{average}\textbf{ }degree $d$ whose AC is at
least $d-2\sqrt{d-1}.$ When $d=3,$ an RSRB$\ \left(  2,6\right)  $ graph is
actually Ramunajan graph with very high probability as $n\rightarrow\infty.$
This is illustrated in Figure \ref{fig:rsrb}(b) where $n\approx2000.$ Out of
1000 randomly chosen such graphs, only \emph{one} had AC less than
$d-2\sqrt{d-1}=0.1716.$

\begin{figure}[tb]
\includegraphics[width=0.97\textwidth]{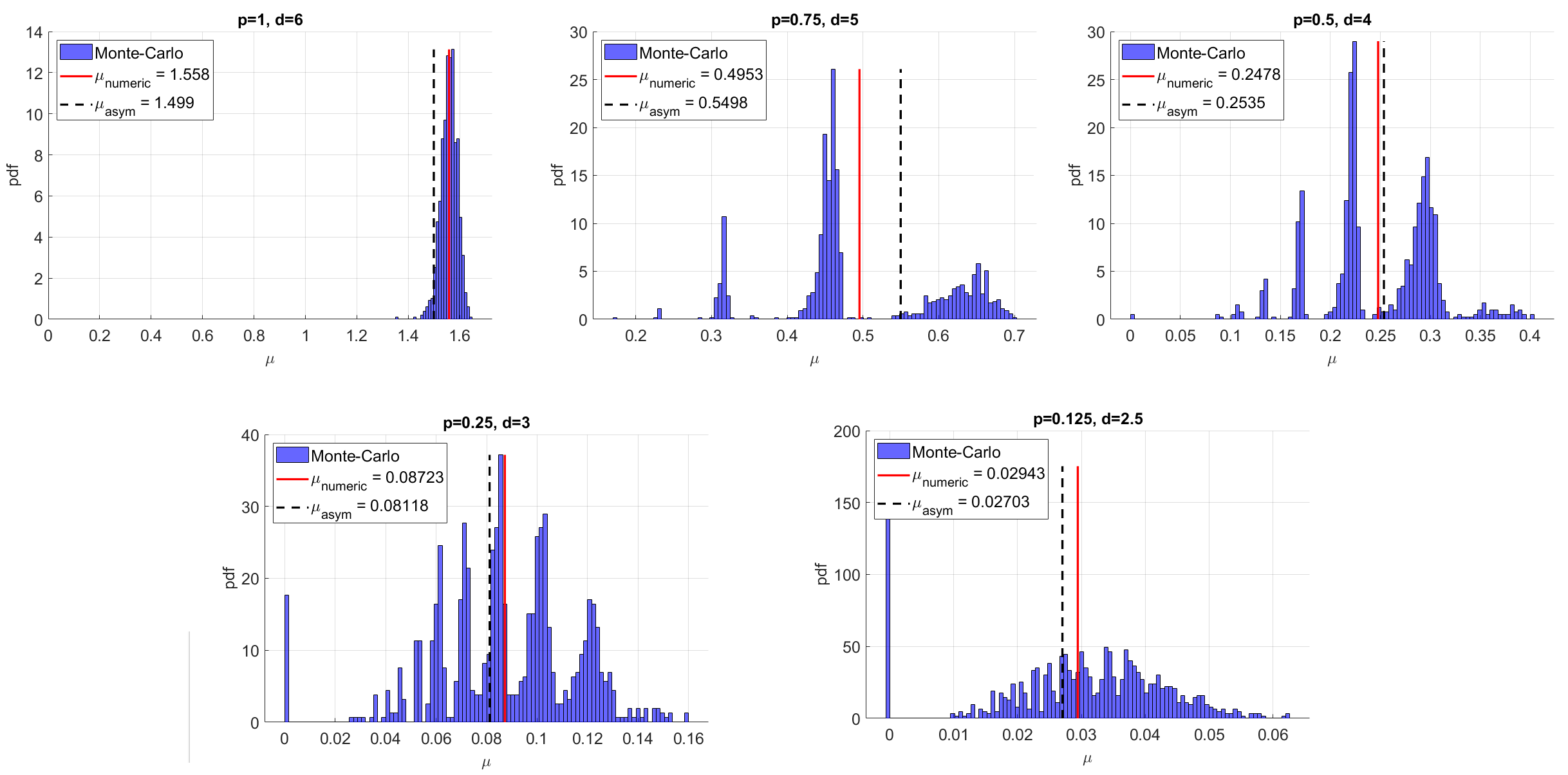}\caption{Distribution of AC
for RSR\ graphs with $\left(  d_{1},d_{2}\right)  =\left(  2,6\right)  $ and
with $p$ as indicated. Average degree $d=2+4p$ is also indicated. We used
$n=500$ nodes, and 1000 simulations. $\mu_{\mbox{numeric}}$ is the average of
the distribution whereas $\mu_{\mbox{asym}}$ is the asymptotics according to
Main Result \ref{main:2}. Note the multi-peaked shape of the distribution, and
the fact that the distribution does not concentrate around the mean when
$d\leq5.$}%
\label{fig:mu26}%
\end{figure}

We relied on careful but formal power series computations. While the results
were validated using numerics, a rigorous analysis of the results in this
paper is an open problem. Related to this, it is an open question to estimate
the accuracy of asymptotics as a function of $n.$ In particular, the accuracy
appears to be significantly degraded for the RSR graphs when $d_{1}$ and
$d_{2}$ far from each other. Take for example the case $\left(  d_{1}%
,d_{2}\right)  =\left(  2,6\right)  $. Intriguingly, the distribution for AC
appears to have multiple peaks and does not necessary concentrate around the
mean, as illustrated in Figure \ref{fig:mu26} (note that this does not appear
to be the case for RSRB graphs as figure \ref{fig:rsrb}(b) illustrates). We
pose this as an open problem.

\textbf{Challenge 1. }\emph{Describe the full distribution of AC, particularly
for RSR graphs. Explain why it can be multi-peaked when }$d_{1}\neq d_{2}%
$.

Another class of graphs having high algebraic connectivity are the complete
bipartite graphs $K_{b,n-b}$, which have AC\ $\mu=b$ (with $b<n/2$). The
average degree of such a graph asymptotes to $d=2b$ as $n\rightarrow\infty,$
so that $\mu\sim d/2.$ Contrast this to $d$-regular random graphs ($\mu\sim
d-2\sqrt{d-1}$). As was noted in \cite{kolokolnikov2015maximizing},
$K_{d/2,n-d/2}$ have higher AC than the $d$-regular graphs (for asymptotically
same number of edges) provided that $d<15$ (since $d/2>d-2\sqrt{d-1}$ in that
case) while the converse is true when $d\geq15.$ While complete bipartite
graphs have a relatively high AC when $d<15$, they are also fragile, in the
sense that removing any single edge causes AC to decrease by one. Random
graphs are more robust with respect to edge deletions as they seem to have
some redundancy built-in.

Consider the limit of large $d_{2}$ for the RSRB\ model. From (\ref{d}) and
(\ref{muRSRB}), one obtains that $d\sim2d_{1}$ and $\mu\sim\frac{d}{2}-1$ as
$d_{2}\rightarrow\infty.$ This is one less than the AC of complete bipartite
graph $K_{d/2,n-d/2}$ (which also has average degree $d$ as $n\rightarrow
\infty)$. More generally, let $\mu_{d_{1},d_{2}}$ be as in\ (\ref{muRSRB}) and
let $\mu_{d}=d-2\sqrt{d-1},$ where $d=2d_{1}d_{2}/(d_{1}+d_{2})$ is the
average degree. It can be shown that $\mu_{d_{1},d_{2}}<\mu_{d}$ for any
$d_{1},d_{2}$ when $d\geq10.$ This discussion suggests the following question.

\textbf{Challenge 2.} \emph{Find a family of random graphs which has a higher
algebraic connectivity than }$d$\emph{-regular random graphs when average
degree }$d\geq10.$ \emph{Explore if more complex degree distribution (e.g.
tri-regular)\ can be better than semi-regular for say, }$d=3.$

There are several notions of graph connectivity, with AC being just one.
Another notion is using the average of reliability polynomial
\cite{brown2014average, perez2018sixty} to compute the effect of edge
deletions on graph connectivity. While the full reliability polynomial
requires exponential time to compute, we performed the following simple
experiment to quickly measure graph robustness. Start with a random 6-regular
graph on 500 vertices (i.e. containing 1500 edges). Then delete edges at
random one by one until the graph becomes disconnected. Over 100 simulations,
it took on average 460 edge deletions until the graph became disconnected (and
in every instance, disconnectivity first occurred when a single vertex lost
all of its edges). We then repeated this experiment with $\left(  d_{1}%
,d_{2}\right)  =\left(  4,12\right)  $ RSRB graph with 500 vertices (which
also contains 1500 edges but has higher AC\ than 6-regular random graph). It
took on average 325 edge deletions until the disconnection was achieved. These
preliminary experiments indicate that $d$-regular graphs are more reliable
with respect to edge deletions than semi-regular bipartite graphs, even though
the latter might have a higher AC. We state this as a conjecture.

\textbf{Conjecture 1.}\emph{ For a given average degree }$d$\emph{, the most
reliable graph is a }$d$\emph{-regular graph.}

Generally speaking, cubic graphs of high girth are good candidates for high-AC
graphs \cite{kolokolnikov2015maximizing}. Since RSRB\ can have higher AC, they
are also natural candidates for searching high-girth graphs. We pose this as a challenge.

\textbf{Challenge 3.}\emph{ For fixed average degree }$d$\emph{ and fixed
number of vertices }$n,$ \emph{find graphs (not necessarily }$d$%
\emph{-regular) with highest possible girth.}

A $d$-regular graph with the smallest possible $n$ for a given girth $g$ is
called a cage. There is an extensive literature for looking for high-girth
graphs; see \cite{exoo2012dynamic, biggs1998constructions} for an overview.
There are powerful methods to do computer searches for high-girth $d$-regular
graphs \cite{mckay1998fast, exoo2011computational}. As an example, consider
$\left(  2,6\right)  $ RSRB\ graphs, whose average degree is $d=3$. Do
$\left(  2,6\right)  $ semi-regular graphs give better cages than cubic graphs
with the same number of vertices? The answer is, it depends.\ Note that that
$\left(  2,6\right)  $ graphs can be obtained from $6$-regular graphs by
inserting a vertex in the middle of each edge; conversely, any $\left(
2,6\right)  $ graph yields a $6$-regular graph by vertex contraction:\ just
delete all 2-degree vertices. According to \cite{o1979smallest,
exoo2012dynamic}, a 6-regular 5-cage graph has 40 vertices. By vertex
insertion, this yields a girth 10 $\left(  2,6\right)  $ semi-regular graph
with 80 vertices. On the other hand, there are three cubic graphs having girth
10 and only 70 vertices \cite{exoo2012dynamic}. So cubic is better for girth 10.

Similarly, a 6-regular cage of girth 6 has 62 vertices \cite{exoo2012dynamic},
which induces a (2,6)\ graph of girth 12 and 124 vertices. On the other hand,
a cubic graph of girth 12 has at least 126 vertices \cite{benson1966minimal,
exoo2012dynamic}. So (2,6) semi-regular graph wins for girth 12.

The techniques in this paper are rather general, and can be used to derive
algebraic connectivity for many other random graph families. As an example,
consider the \textquotedblleft small-world\textquotedblright-type graph
illustrated in Figure \ref{fig:graphs}(c). Start with a cycle of $n$ vertices.
Then connect all odd-numbered vertices to each other at random. See Appendix A
for Matlab code that generates such a graph. The resulting graph has an
average degree of 2.5. To compute its AC, as in the derivation of
(\ref{muRSRB}), we add a self-loop to even-numbered vertices so that each
vertex has degree 3. We leave it as an excercise to the reader that the
average walk-counting sequence $\phi_{s}$ has the generating function which
satisfies%
\begin{align}
A &  =1+x^{2}(A\hat{A}+AB)\nonumber\\
\hat{A} &  =1+2x^{2}\hat{A}B\nonumber\\
B &  =1+xB+x^{2}AB\nonumber\\
\phi_{o} &  =1+x^{2}\phi_{o}\left(  \hat{A}+2B\right)  \\
\phi_{e} &  =1+x\phi_{e}+2x^{2}\phi_{e}A\nonumber\\
\phi &  =\frac{1}{2}\phi_{o}+\frac{1}{2}\phi_{e}\nonumber
\end{align}
Looking at the radius of convergence of the generating functions as before, we
find that in the limit $n\rightarrow\infty$ its AC asymptotes to the smallest
root of the polynomial%
\[
\mu^{9}-24\mu^{8}+249\mu^{7}-1454\mu^{6}+5184\mu^{5}-11400\mu^{4}+14848\mu
^{3}-10368\mu^{2}+3108\mu-136=0,
\]
explicitly given by $\mu=0.0521926.$ This is validated using Monte Carlo
simulations:\ the average AC\ of 1000 random such graphs with $n=1000$
vertices is 0.0557 (std=0.0029), the difference of about 7\%. This is
significantly higher than a RSR\ graph with $(p,d_{1},d_{2})=\left(
0.5,2,3\right)  $, which also has average degree of 2.5, but whose AC\ is
$\mu\sim0.044241.$

\section*{Appendix A: Matlab and Maple code}

\vskip 2cm

\begin{minipage}[h]{0.45\textwidth}
\textbf{RSRB model (Figure \ref{fig:graphs}(a)):}\newline
\newline
\texttt{d1=2; d2=3; n1=30; n2=20;}\newline%
\texttt{bag1=mod([0:n1*d1-1], n1)+1;}\newline\texttt{bag2=mod([0:n2*d2-1],
n2)+1+n1;}\newline\texttt{bag2=bag2(randperm(numel(bag2)));}\newline%
\texttt{G=graph(bag1, bag2);}\newline\texttt{plot(G);}
\end{minipage}

\vskip 2cm \begin{minipage}[h]{0.45\textwidth}
$~$\newline
\textbf{RSR model (Figure \ref{fig:graphs}(b)):}\newline
\newline
\texttt{p=0.4; d1=2; d2=3; n=50;}\newline%
\texttt{n1=(1-p)*n; n2=p*n;}\newline\texttt{v1=mod([0:n1*d1-1], n1)+1;}%
\newline\texttt{v2=mod([0:n2*d2-1], n2)+1+n1;}\newline\texttt{bag=[v1,
v2];}\newline\texttt{bag=bag(randperm(numel(bag)));}\newline%
\texttt{G=graph(bag(1:end/2), bag(end/2+1:end));}\newline\texttt{plot(G);}
\end{minipage}

\vskip 2cm \begin{minipage}[h]{0.45\textwidth}
$~$\newline
\textbf{Small-world-like model (Figure \ref{fig:graphs}(c)):}\newline
\texttt{
\\
n=52;\\
v1=[1:n];\\
v2=[2:n,1];\\
bag=[1:2:n];\\
bag=bag(randperm(numel(bag)));\\
v1=[v1,bag(1:end/2)];\\
v2=[v2,bag(end/2+1:end)];\\
G=graph(v1,v2);\\
t=(1:n)/n*2*pi;\\
plot(G, 'XData', cos(t), 'YData', sin(t));
}
\end{minipage}

\bigskip

\begin{absolutelynopagebreak}
\textbf{Maple code to compute AC\ of RSR\ graphs:}\newline\newline\includegraphics[width=0.8\textwidth]{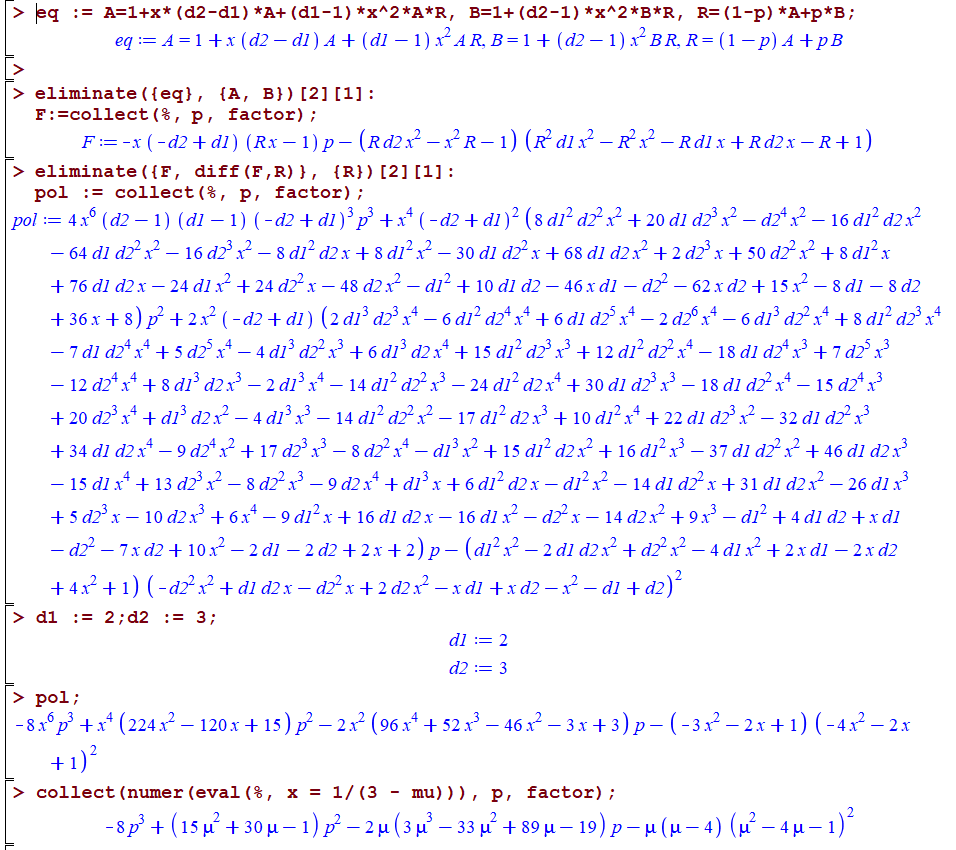}
\end{absolutelynopagebreak}

\bibliographystyle{elsarticle-num}
\bibliography{semireg}

\end{document}